\documentclass[11pt,twoside]{article}
\usepackage{comment}
\usepackage{amsmath, amssymb, amsfonts, amstext, amsthm, textcomp, enumerate}
\usepackage[mathscr]{euscript}
\usepackage{float}
\usepackage{epsfig,epsf,epstopdf}
\usepackage{relsize}
\usepackage[pdf]{pstricks}
\usepackage{float}
\usepackage{booktabs}
\usepackage{booktabs}
\usepackage{mathtools}
\usepackage{tabularx}
\usepackage[left=20mm,top=0.5in,bottom=15mm]{geometry}
\usepackage{graphicx}
\usepackage{caption}
\usepackage{epstopdf}
\usepackage{longtable}
\usepackage[utf8]{inputenc}
\usepackage{color}
\usepackage{hyperref}
\usepackage{graphicx}
\usepackage{dcolumn}
\usepackage{bm}
\usepackage{epstopdf}
\usepackage[english]{babel}
\usepackage{subfigure}
\usepackage{color}

\usepackage{ulem}
\usepackage{tikz}
\usepackage{bbm}

\newtheorem{thm}{Theorem}[section]

\newtheorem{prob}[thm]{Problem}
\newtheorem{lem}[thm]{Lemma}
\newtheorem{defn}[thm]{Definition}

\newtheorem{rem}[thm]{Remark}
\newtheorem{ex}[thm]{Example}

\newcommand{\mnr}{\mathbf{M}_n\,(\mathbb{R})}

\texorpdfstring

\newfont{\bb}{msbm10}

\def\p{{\large\textbf{P} }}

\title{Open Problems on the Null Space and Image Space of the Adjacency Matrix}
\author{K. Kranthi Priya \thanks{Department of Mathematics, Indian Institute of Technology Madras, Chennai 600036, India (ma21d010@smail.iitm.ac.in, kcskumar@iitm.ac.in).}  $\cdot$  K. C. Sivakumar $^\dagger$ }
\begin{document}
\maketitle
\begin{abstract}
Akbari, Cameron, and Khosrovshahi proposed two questions concerning the null space and the image space of graphs. Although these problems have been studied in several cases, they have not yet been fully resolved. In this work, we present constructions of graphs that satisfy both the problems. 

\end{abstract}


\newpage
\section{Introduction}

Let $G$ be a graph with vertex set 
$
V(G) = \{v_1, v_2, \ldots, v_n\},$
and edge set
\[
E(G) = \{(v_i, v_j) : v_i \sim v_j\},
\]
where the notation $v_i \sim v_j$ indicates that the vertices $v_i$ and $v_j$ are {\it adjacent} in $G$, i.e., there is an edge between $v_i$ and $v_j$. The order and size of the graph $G$ are $|V(G)|$ and $|E(G)|$, respectively.

For a graph $G$, the adjacency matrix is a square matrix of order $n$ such that the $(i,j)$ entry of $A(G)$ is $1$ if the vertices $v_i$ and $v_j$ are adjacent and $0$ otherwise, with respect to an arbitrary labeling of the vertices of $G.$ We often denote $A(G)$ simply by $A.$ The {\it nullity} $n(G)$ and {\it rank} $r(G)$ of a graph $G$ are the nullity and rank of the associated matrix $A(G),$ respectively.

Characterizing graphs $G$ for which $n(G) >0$, is of great
interest in chemistry \cite{LH50}. In particular, for bipartite graphs representing  alternate hydrocarbons, the condition $n(G)>0$ indicates molecular 
instability. This problem, originally posed by Collatz and Sinogowitz \cite{CS57}, remains open in general, with complete results known only for certain classes of graphs, including trees and some 
bipartite graphs \cite{FDS05,GS01}.

The problem of bounding the order of a graph in terms of the
rank was first investigated by Kotlov and Lovás \cite{KL96}. A graph is said to be {\it reduced} if no two
vertices have the same set of neighbours. Let $m(r)$ denote the maximum possible order of
a reduced graph of rank $r$. In \cite{KL96}, authors proved that there exists a constant $c$ such that every reduced graph of rank $r$ has order at most $c2^{r/2}$, Moreover, for
every $r\ge 2$, they constructed a graph of rank $r$ and order $n(r)$, where
\begin{equation}
    n(r)=\left\{\begin{array}{cl} 2^{\frac{r+2}{2}}-2,
         & r \quad\text{is even} \\[2ex]
         5\cdot{}2^{\frac{r-3}{2}}-2,
         &  r\quad \text{is odd}.
    \end{array}
    \right. 
\end{equation}
Akbari et al. \cite{ACK} conjectured that, in fact, $m(r) = n(r)$. The following problems were also proposed by Akbari et al.

\begin{prob}\label{conj1}\cite[Question~1]{ACK}
    For which graphs there is a $\{-1,0,1\}$-basis for the null space?
\end{prob}
 
\begin{prob}\label{conj2}\cite[Question~2]{ACK}
 Given any graph $G$ with at
least one edge, there exists a nonzero $(0, 1)$-vector in the row space of its adjacency matrix $A$ (over the real
numbers), which does not occur as a row of $A.$
\end{prob}

The problem \ref{conj2} is also known as ACK conjecture. The motivation for studying the ACK conjecture stems from the observation that, if the conjecture holds for all graphs, then Theorem 3 in \cite{BBDG97} implies that the function $m(r)$ is increasing. This result, in turn, provides further support for the conjectured equality $m(r)=n(r)$. In addition to the ACK conjecture, if we find a $(0,1)$ vector $y$ such that $A(G)x=y$ and $x^Ty=0.$ Then, using Theorem 3 in \cite{BBDG97} proves that the rank of the new graph, which is obtained by adding one vertex to $G$, joins a
suitable set of existing vertices so that the rank increases by only one. This establishes that $m(r)$ is strictly monotonic, which is a conjecture in \cite{ACK}.
   
 Both problems haven't been resolved completely.
The authors of~\cite{AAGK05} and~\cite{SS05} independently proved that problem \ref{conj1} holds for forests. Moreover, in~\cite{JMPS18}, an algorithm was developed to determine a sparsest basis of the null space of any given forest, which is additionally a $\{-1,0,1\}$-basis. Interestingly, the problem does not hold for unicyclic graphs; consequently, T. Sander et al. \cite{SS07} characterized all unicyclic graphs satisfying problem \ref{conj1}. They also established problem for cographs \cite{S03}. Nevertheless, the problem remains open.

 In \cite{B25}, the author proves that every graph with diameter at least $4$ satisfies problem \ref{conj2}, and also addresses few cases for graphs of diameter $2$ and $3$. Furthermore, \cite{SKBD25} identifies a class of graphs that may contain potential counter examples. Interestingly, in \cite{AK26}, the authors construct families of graphs belonging to this potential counterexample class that nevertheless satisfy problem \ref{conj2}.

In this paper, we present families of graphs that satisfy both problems~\ref{conj1} and~\ref{conj2}, while also belonging to the potential counterexample class. The paper is organised as follows. Section~2 introduces the necessary definitions and results in the literature. In Section~3, we present constructions of core graph classes that satisfy both Problems~\ref{conj1} and~\ref{conj2}. Moreover, the constructions in this section belong to the potential counter-example class. Finally, in the last section, we show that every nut graph of diameter 
$3$ satisfies the problem~\ref{conj2}.

\section{Basic Definitions and Notations}

Let $\mnr$ be set of all real square matrices of order $n$. For $A \in \mnr,$ we denote by $N(A),$ the {\it null space} of $A$ and use $R(A),$ to denote the {\it range space} of $A$. The dimension of the subspace $R(A)$ is called the {\it rank of $A$} and is denoted by $\rm r(A)$. The notation $A(\alpha)$ represents the submatrix of $A$ formed by removing rows and columns indexed by $\alpha$, where $\alpha$ is a subset of $\{1, 2,\ldots, n\}.$ 
\begin{defn}\label{detblockfor}
    Let $A\in \mnr$ and be partitioned as
    $$A=\begin{pmatrix}
        A_{11} & A_{12}\\
        A_{21} & A_{22}
    \end{pmatrix},$$
    where $A_{11}$ and $A_{22}$ are square matrices. If $A_{11}$ is invertible then the {\it Schur complement} of $A_{11}$ in $A$ is defined to be the matrix $A_{22}-A_{21}A_{11}^{-1}A_{12}.$ Similarly, if $A_{22}$ is nonsingular then the Schur complement of $A_{22}$ in $A$ is defined to be the matrix $A_{11}-A_{12}A_{22}^{-1}A_{21}.$  
    \end{defn}
It is easy to obtain the following:
\begin{equation}\label{blockexp}
    \begin{pmatrix}
        I & {\bf 0}\\
        -A_{21}A_{11}^{-1} & I
\end{pmatrix}\begin{pmatrix}
        A_{11} & A_{12}\\
        A_{21} & A_{22}
\end{pmatrix}\begin{pmatrix}
        I & -A_{11}^{-1}A_{12}\\
        {\bf 0} & I
    \end{pmatrix}=\begin{pmatrix}
        A_{11} & {\bf 0}\\
        {\bf 0} & A_{22}-A_{21}A_{11}^{-1}A_{12}
    \end{pmatrix}.
\end{equation}

Using the above fact, we conclude \begin{equation}\label{detfor}
\det(A)=\det(A_{11})\det(A_{22}-A_{21}A_{11}^{-1}A_{12}).
\end{equation}

 Let $G$ be a graph with vertex set $V(G)$. For a vertex $v\in V(G)$, we denote $N(v)$ as the set of vertices that are adjacent to vertex $v$ in $G,$ usually called the {\it Neighbourhood} of $v$ in $G$. The {\it degree of a vertex} $v$, denoted by $d(v)$, is the cardinality of the set $N(v)$. The vertex $v$ is said to be a {\it pendant} vertex if $d(v)=1.$ A {\it path} of length $k$, denoted by $P_{k+1}$ between two vertices $v_1,v_{k+1}$ is an alternating sequence of vertices and edges $v_1,e_1,v_2,e_2,\ldots,e_{k},v_{k+1},$ where each vertex $v_i$ is distinct and $e_i$ is an edge $(v_i,v_{i+1})$. 
A graph is said to be {\it connected} if there is a path between any pair of vertices in $G.$ The distance
between two vertices $u$ and $v$ in a connected graph $G$ is the length of the shortest path between them, and it is denoted by $d(u, v).$ Clearly, if $u$ and $v$ are adjacent, then
$d(u, v) = 1.$ For a graph $G$, its diameter is defined as $diam(G) = \max_{u,v\in V} d(u, v).$  Throughout this paper, we restrict to {\it simple} graphs, which are graphs with no multiple edges or loops. A graph is said to be a {\it d-regular graph} if every vertex has the same degree $d.$

A {\it path graph} $P_n$ is a connected graph on $n$ vertices that has exactly two pendant vertices, and the degree of the other vertices is $2.$  
The determinant of the adjacency matrix of a path graph $P_n$ is obtained as follows in \cite{BDH98}:
\begin{eqnarray}\label{detpath}
 \det(A(P_n))=  \begin{dcases}
     1 & n\cong 0(mod~4)\\
     -1 & n\cong 2(mod~4)\\
     0 & \text{n is odd.}
 \end{dcases} 
\end{eqnarray}
The adjacency matrix of a path is nonsingular if and only if the path has an even number of vertices.

A {\it subgraph} of a graph 
$G$ is a graph $\Gamma$ whose vertex set satisfies 
$V(\Gamma)\subseteq V(G)$ and whose edge set satisfies $E(\Gamma)\subseteq E(G)$. An {\it induced subgraph} of $G$ is a subgraph $\Gamma$ in which every edge 
$(v_i,v_j)$ with 
$v_i,v_j\in V(\Gamma)$ that appears in 
$G$ also appears in $\Gamma$.

A vertex $v_i$ in a graph $G$ is called a {\it core vertex} (CV) if there exists $x=(x_i)\in N(A(G))$ such that $x_i,$ is nonzero.  A graph 
$G$ is called a {\it core} graph if every vertex of 
$G$ is a core vertex. 

The graph $G+v$ is obtained from the graph $G$ by adding a vertex $v,$ that is isolated or joined by one or more edges to $G.$ There is another definition of a core vertex depending on the change in nullity produced with this operation. The vertex $v$ is a core vertex (also referred to as a {\it downer} vertex in the literature) in $G+v$ if $n(G)=n(G+v)-1.$ 
 
An eigenvalue of a matrix is {\it main}, if there exists an associated eigenvector $x$ not orthogonal to the all-one vector ${\bf e}.$


The problem \ref{conj2} has been proven for all graphs of diameter larger than $4.$

\begin{thm}\cite[Theorem 3.1]{B25}
Let $G$ be a simple connected graph with $diam(G)\ge 4.$ Then there
exists a nonzero $(0, 1)$-vector in the row space of $A(G)$ (over the real numbers) which
does not occur as a row of $A(G).$
\end{thm}

The following results prove that certain classes of graphs with diameters 
$2$ and 
$3$ satisfy problem \ref{conj2}.

\begin{thm}\cite[Theorem 3.2]{B25}
    Let $G$ be a simple non-complete connected graph such that $diam(G) \le 3$
and $r(G) \le 5$. Then there exists a non-zero $(0, 1)$-vector in the row space of $A(G)$
(over the real numbers) which does not occur as a row of $A(G)$.
\end{thm}

\begin{thm}\cite[Theorem 3.3]{B25}\label{e2n-5}
    Let $G$ be a diameter-2 graph of order $n$ with no dominating vertex.
Suppose $|E(G)| = 2n-5.$ Then there exists a non-zero $(0, 1)$-vector in the row space of
$A(G)$ (over the real numbers), which does not occur as a row of $A(G)$.
\end{thm}

Let ${\cal C}$ denote the class of all graphs $G$ which have the potential to be counter-examples to the problem \ref{conj2}. Recently, the following necessary conditions for such graphs were proved. 

\begin{thm}\label{sciriha}
\cite[Theorem 5.1]{SKBD25}
Let $G \in {\cal C}$. We then have the following:
\begin{enumerate}
\item $G$ is a core graph;
\item $0$ is a main eigenvalue of $G$;
\item every vertex of $G$ lies on a triangle;
\item for every vertex $u \in V(G)$, each vertex in $N(u)$ forms a triangle containing $u$;
\item $G$ is not regular;
\item $G$ must be connected;
\item $G$ is not bipartite;
\item the diameter of $G$ is $2$ or $3$.    
\end{enumerate} 
\end{thm}

\section{Construction of core graphs}

\subsection{Construction~I}
\begin{defn}\label{Corewheel1}
Let $n$ be a positive integer. We introduce a family of graphs $\mathscr{C}_1$ whose vertex set $V$  is partitioned into three pairwise disjoint subsets $V_0, V_1$ and $V_2$. Specifically, we set 
\begin{center}
 $V_0=\{u_1,\ldots, u_{4k}\},~V_1=\{u_{4k+1}\},$ for some $k.$    
\end{center}

The edge set \( E \) of \( \mathscr{C}_1\) is defined as follows:
\begin{enumerate}
    \item The vertices in $V_0 \cup V_1$ induces a wheel graph, where $u_{4k+1}\in V_1$ serves as the centre. 
    \item The adjacency relation for vertices in $V_2$ satisfies the following:
    \begin{itemize}
        \item [$(i)$] $\{N(v): v\in V_2 \}$ form a partition of $V_0.$
        \item [$(ii)$] For each $v\in V_2,$ $|N(v)|=4t$ for some $t.$
    \item[$(iii)$]  For each $v\in V_2$, if $u_i,u_j\in N(v)$ with $i<j<4k+1$, then $u_{\ell}\in N(v)$ for $i\le \ell  \le j.$ Additionally, both $u_1,u_{4k}$ together does not belong to $N(v)$ for any $v\in V_2$, except the case $|V_2|=1$.
\end{itemize}
 \end{enumerate}
\end{defn}

 \begin{figure}[H]
\begin{center}
\begin{tikzpicture}[scale=1]
\draw  (5.,3.)-- (6.,5.);
\draw  (6.,5.)-- (7.,3.7);
\draw  (7.,3.7)-- (5.,3.);
\draw  (5.,3.)-- (7.,2.3);
\draw  (7.,2.3)-- (6.,1.);
\draw  (6.,1.)-- (5.,3.);
\draw  (5.,3.)-- (4.,1.);
\draw  (4.,1.)-- (3.,2.4);
\draw  (3.,2.4)-- (5.,3.);
\draw  (5.,3.)-- (3.,3.7);
\draw  (3.,3.7)-- (4.,5.);
\draw  (4.,5.)-- (5.,3.);
\draw  (6.,5.) -- (4.,5.);
\draw  (6.,1.)--(4.,1.);
 \draw  (3.,3.7)--(3.,2.4);
 \draw (7.,3.7)--(7.,2.3) ;
 \draw  (.5,3.)--(3.,3.7);
 \draw (.5,3.)--(4.,5.);
 \draw (.5,3.)--(4.,1.);
 \draw (.5,3.)--(3.,2.4);
 \draw (10.,3.)--(6.,5.);
  \draw (10.,3.)--(7.,3.7);
   \draw (10.,3.)--(7.,2.3); \draw (10.,3.)--(6.,1.);
\begin{scriptsize}
\fill (10.,3.) circle (2.5pt);
\draw (10.35,3.) node {$11$};
\fill (.5,3.) circle (2.5pt);
\draw (.2,3.)
node {$10$};
\fill (5.,3.) circle (2.5pt);
\draw (5.,3.45) node {$9$};
\fill (6.,5.) circle (2.5pt);
\draw (6.25,5.25) node {$1$};
\fill (7.,3.7) circle (2.5pt);
\draw (7.2,3.9) node {$2$};
\fill (7.,2.3) circle (2.5pt);
\draw (7.2,2.1) node {$3$};
\fill (6.,1.) circle (2.5pt);
\draw (6.25,0.8) node {$4$};
\fill (4.,1.) circle (2.5pt);
\draw (3.75,0.7) node {$5$};
\fill (3.,2.4) circle (2.5pt);
\draw (2.6,2.) node {$6$};
\fill (4.,5.) circle (2.5pt);
\draw (3.9,5.3) node {$8$};
\fill (3.,3.7) circle (2.5pt);
\draw (2.7,4.) node {$7$};
\end{scriptsize}
\end{tikzpicture}
\end{center}
\caption{}
\label{fig:7}
\end{figure}

\begin{lem}\label{C1diam}
        Let $G\in$ \( \mathscr{C}_1\). If $|V_2|\le 3$ then $diam(G)=3.$ 
\end{lem}

\begin{proof}
Let $u,v\in V$ then it is enough to show that $d(u,v)\le 3.$ \\
     Suppose $u,v$ are adjacent then there is nothing to prove. Let us assume $u, v$ are not adjacent. \\{\bf Case~1}: For $u,v\in V_0\cup \{u_{4k+1}\}.$\\
     Since $V_0\cup \{u_{4k+1}\}$ is a wheel graph, $d(u,v)\le 2.$ \\
    {\bf Case~2}: For $u\in V_0\cup \{u_{4k+1}\},~v\in V_2.$\\
     Suppose $u=u_{4k+1}$ and $v\in V_2$ then using $2(i)$ of the definition \ref{Corewheel1}, there exists $u_{t}\in V_0$ such that $u_t\in N(v)$. 
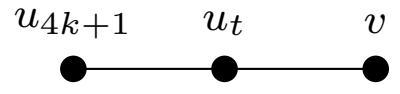
     \begin{figure}[H]
    \begin{center}
\begin{tikzpicture}[scale=1]
   \draw (7.,3.)--(6.,3.);
   \draw (7.,3.)--(8.,3.);
   
\begin{scriptsize}

\fill (6.,3.) circle (2.5pt);
\draw (6.,3.3) node {$u_{4k+1}$};
\fill (7.,3.) circle (2.5pt);
\draw (7.,3.3) node {$u_t$};
\fill (8.,3.) circle (2.5pt);
\draw (8.,3.3) node {$v$};

\end{scriptsize}
\end{tikzpicture}
\end{center}
\end{figure}
     Thus, we can conclude that $d(u_{4k+1},v)=2.$ On the other hand, if $u\in V_0$ and $v\in V_2$. We know that $u$ and $u_{4k+1}$ are adjacent. From the above argument we have $d(u_{4k+1},v)=2$. Therefore, we conclude $d(u,v)\le 3.$  \\{ \bf Case~3}: For $u,v\in V_2$.\\
     We have $|V_2|=k\le 3.$ Suppose $i_0,j_0~(i_0<j_0)$ be the least indices among all the neighbours of $u,v$ respectively in $V_0.$ If $u_{j_0-1}\in N(u)$ then we conclude that there exists a path from $u,v$ passing through $u_{j_0-1},u_{j_0}$, which is of length $3$. i.e., 
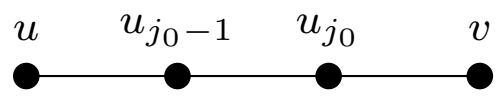
\begin{figure}[H]
\begin{center}
\begin{tikzpicture}[scale=1]
\draw (5.,3.)--(6.,3.);
   \draw (7.,3.)--(6.,3.);
   \draw (7.,3.)--(8.,3.);
   
\begin{scriptsize}
\fill (5.,3.) circle (2.5pt);
\draw (5.,3.3) node {$u$};
\fill (6.,3.) circle (2.5pt);
\draw (6.,3.3) node {$u_{j_0-1}$};
\fill (7.,3.) circle (2.5pt);
\draw (7.,3.3) node {$u_{j_0}$};
\fill (8.,3.) circle (2.5pt);
\draw (8.,3.3) node {$v$};

\end{scriptsize}
\end{tikzpicture}
\end{center}
\end{figure} 
If $u_{j_0-1}\notin N(u)$, then we have two cases to deal with. Suppose $i_0>1$ then there exists a path from $u,v$ passing through $u_{i_0-1},u_{i_0},$ that is of length $3.$
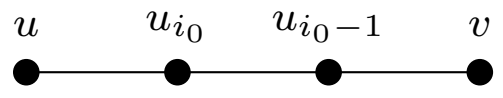
\begin{figure}[H]
\begin{center}
\begin{tikzpicture}[scale=1]
\draw (5.,3.)--(6.,3.);
   \draw (7.,3.)--(6.,3.);
   \draw (7.,3.)--(8.,3.);
   
\begin{scriptsize}
\fill (5.,3.) circle (2.5pt);
\draw (5.,3.3) node {$u$};
\fill (6.,3.) circle (2.5pt);
\draw (6.,3.3) node {$u_{i_0}$};
\fill (7.,3.) circle (2.5pt);
\draw (7.,3.3) node {$u_{i_0-1}$};
\fill (8.,3.) circle (2.5pt);
\draw (8.,3.3) node {$v$};
\end{scriptsize}
\end{tikzpicture}
\end{center}
\end{figure}  On the other hand, let us say $i_0=1$ then there exists a path from $u,v$ passing through $u_{4k},u_{1},$ that is of length $3.$ Thus $d(u,v)=3.$
\end{proof}

The following lemmas are related to the image space of the adjacency matrix of a path graph. These results are used to prove Theorem \ref{C1core}.
\begin{lem}\label{lem1}
    Let $A$ be the adjacency matrix of the path graph $P_{2h}.$
         If $Ay=\begin{pmatrix}{\bf 0}_{m}
    \\{\bf 1}_{4s} \\ {\bf 0} \end{pmatrix}$ then, $y=\begin{pmatrix}
        {\bf 0}_m\\
        \alpha_s\\
        {\bf 0}
    \end{pmatrix}$
    where, $\alpha_s=(\alpha,\ldots,\alpha)^T$ and $ \alpha=(0,1,1,0)^T$.
    
\end{lem}

\begin{proof}

    It is enough to prove this claim for $s=1,$ due to fact that $A$ is linear over the vectors. Thus, $\alpha_s=\alpha$ and we obtain, $$A\begin{pmatrix}
        {\bf 0}_m\\
        \alpha\\
        {\bf 0}
    \end{pmatrix}=Ae^{m+2}+Ae^{m+3}=\begin{pmatrix}
        {\bf 0}_{m} \\
        1\\
        0\\
        1\\
        0\\
        {\bf 0}
    \end{pmatrix}+\begin{pmatrix}
        {\bf 0}_{m}\\
        0\\
        1\\
        0\\
        1\\
        {\bf 0}
    \end{pmatrix}=\begin{pmatrix}
        {\bf 0}_m \\
        1\\
        1\\
        1\\
        1\\
        {\bf 0}
    \end{pmatrix}=\begin{pmatrix}{\bf 0}_{m}
    \\{\bf 1}_{4} \\ {\bf 0} \end{pmatrix}.$$
      \end{proof}

     \begin{lem}\label{pathodd}
     Let $A$ be the adjacency matrix of the path graph $P_{2h}$. If $h$ is odd, then the following statements hold:
    \begin{enumerate}
        \item  If $Ay=\begin{pmatrix}
        1\\
        1\\
        {\bf 0}
    \end{pmatrix}$ 
    then $y=  \begin{pmatrix}
        1\\
        1\\\zeta_{\frac{h-1}{2}}
   \end{pmatrix}$      where, $\zeta_s=(\zeta,\ldots,\zeta)^T$ and $    \zeta=(0,-1,0,1)^T.$
    \item  If  $Ay=\begin{pmatrix}{\bf 1}_{4t+2} \\ {\bf 0} \end{pmatrix}$ then
    $y=\begin{pmatrix}
        1\\
        1\\
        \beta_{t}\\
      \zeta_{\frac{h-(2t+1)}{2}}
    \end{pmatrix}$
    where, $\beta_s=(\beta,\ldots,\beta)^T$ and \\$\beta=(0,0,1,1).$
     
    \end{enumerate}
    \end{lem}
    \begin{proof}
       \begin{enumerate}
    \item Here, $h$ is odd, so $|V(G)|=2h=2(2\ell+1)=4\ell+2$ for some $\ell.$
     Let us prove the result using the induction hypothesis on $\ell.$
\\For $\ell=1$, $|V(G)|=6$. We have, $$A(P_6)=\begin{pmatrix}
  0 & 1 & 0 & 0 & 0 & 0\\
         1 & 0 & 1 & 0& 0 &0\\
         0 & 1 & 0 & 1& 0 &0\\
         0 & 0 & 1 & 0 & 1 &0\\
         0 & 0 & 0 & 1 & 0 & 1 \\
         0 & 0 & 0 & 0  & 1 & 0
\end{pmatrix}$$
     It is easy to verify that $$A(P_6)\begin{pmatrix}
         ~~1\\
         ~~1\\
         ~~0\\
         -1\\
         ~~0\\
         ~~1
     \end{pmatrix}=\begin{pmatrix}
         1 \\
         1\\
         0\\0\\0\\0
     \end{pmatrix}.$$ Let us assume that the claim holds for a path graph of length $4(\ell-1)+2=4\ell-2.$ Thus, there exists a vector $d$ such that the following holds: i.e., $$A(P_{4\ell-2})d=A(P_{4\ell-2})\begin{pmatrix}
         1\\
         1\\
\zeta_{\frac{2\ell-2}{2}}
     \end{pmatrix}=e^1+e^2.$$ For any path graph of length $4\ell+2$, we have $$A(P_{4\ell+2})=\left[
\begin{tabular}[c]{c|cccc}
$A(P_{4\ell-2})$ & $e^{4\ell-2}$ & ${\bf 0}$ & ${\bf 0}$& ${\bf 0}$\\
\hline
$(e^{4\ell-2})^T$ & $0$ & $1$ &  $0$ & $0$\\
${\bf 0}^T$ & $1$ & $0$ & $1$ &$0$\\
${\bf 0}^T$ & $0$ & $1$ & $0$ &$1$\\${\bf 0}^T$ & $0$ & $0$ & $1$ &$0$\\
\end{tabular}
\right].$$
  Therefore, \begin{eqnarray*}
    A(P_{4\ell+2})\begin{pmatrix}
      1\\
         1\\
        \zeta_{\frac{(2\ell+1)-1}{2}}
  \end{pmatrix} =A(P_{4\ell+2})\begin{pmatrix}
      1\\
         1\\
         \zeta_{\frac{2\ell-2}{2}}\\
         \zeta
  \end{pmatrix} &=&\begin{pmatrix}
      A(P_{4\ell-2})d\\
     ( e^{4\ell-2})^Td-1\\
     0\\
      0\\
      0
  \end{pmatrix} \end{eqnarray*}
  We can observe that the last entry of $d$ is $1$. Thus, $( e^{4\ell-2})^Td-1=0$. Using the fact that $A(P_{4\ell-2})d=e^1+e^2$. Therefore, the hypothesis holds for the general case.

     \item 
    We also have $Ay=\begin{pmatrix}
        1\\
        1\\0\\
        {\bf 0}
    \end{pmatrix}+\begin{pmatrix} 0\\
    0
    \\{\bf 1}_{4t} \\ {\bf 0} \end{pmatrix}.$ 
    Therefore, using the fact $A$ is invertible and $(1),(2)$, we obtain the following:$$y=A^{-1}\begin{pmatrix}
        1\\
        1\\0\\
        {\bf 0}
    \end{pmatrix}+\begin{pmatrix} 0\\
    0
    \\ \alpha_{t} \\ {\bf 0} \end{pmatrix}=\begin{pmatrix}
        1\\
        1\\\zeta_{t}\\
        \zeta_{\frac{h-2t-1}{2}}
   \end{pmatrix}+\begin{pmatrix} 0\\
    0
    \\ \alpha_{t} \\ {\bf 0} \end{pmatrix}=\begin{pmatrix}
        1\\1\\\alpha_{t}+\zeta_{t}\\
\zeta_{\frac{h-2t-1}{2}}
    \end{pmatrix},$$
    where $f_{t}=(f,\ldots,f)$ and $\zeta=(0,-1,0,1),~\alpha=(0,1,1,0)$. Then $\alpha+\zeta=(0,0,1,1),$ establishes the claim.
    \end{enumerate} 
    \end{proof}

\begin{lem}\label{cyclenull}
    Let $G$ be a cycle $C_{4k}$ defined on the vertex set $\{u_1,\ldots,u_{4k}\}$ for some $k$. Then $$\{(1,0,-1,0,1,0,\ldots,-1,0)^T,(0,1,0,-1,0,\ldots,0,-1)^T\}$$ form a basis for $N(A)$. 
\end{lem}

\begin{proof}
We have $$A(C_{4k})=\begin{pmatrix}
0  & 1 & ~ & ~ & 1\\
1 & 0 & 1 & ~ & ~ \\
~ & 1 & \ddots &\ddots& ~ \\~ & ~ & \ddots & 0 & 1 \\
1 & ~ & ~ & 1 & 0 \\
\end{pmatrix}$$ 
Note that the entries that are not filled are zero. Suppose $x=(x_i)\in N(A)$ then we have 
  $x_i+x_j=0$ for $|i-j|=2$, $x_1+x_{4k-1}=0$ and $x_2+x_{4k}=0.$
Consequently, $$N(A)=\{(x_i): x_{2t}=(-1)^{t+1}x_2, ~ x_{2t+1}=(-1)^{t+2}x_1\}.$$ This establishes the result.
\end{proof}

\begin{thm}\label{C1core}
      Every graph $G\in$ \( \mathscr{C}_1\) is a core graph.
\end{thm}
  \begin{proof} 
  The sets $V_0,V_1$ are defined as in definition \ref{Corewheel1}. Using Lemma \ref{C1diam}, we restrict $|V_2|\le 3.$ We denote $$V_2=\{u_{4k+2},u_{4k+3},u_{4k+4}\}.$$
  We have $\{N(v): v\in V_2 \}$ form a partition of $V_0.$ Thus, 
we conclude that \begin{eqnarray}\label{equ1v2}
    -a^{4k+1}+\sum_{\{i~:~u_i\in V_2\}}a^{i}=0.
\end{eqnarray} Thus, $z=(z_i)\in N(A(G))$, where, $z_i=\begin{dcases}
    -1 & i={4k+1},\\
    1 & 4k+2\le i\le 4k+4 ,\\
    0 & \text{Otherwise}.
    \end{dcases}$

   We have $$A(G)=\begin{pmatrix}
        A(C_{4k}) & {\bf 1} & {\bf w}_1 & {\bf w}_2 & {\bf w}_3\\
       { \bf 1}^T & 0 & 0 & 0 & 0\\
        {\bf w}^T_1 & 0 & 0 & 0 & 0\\
        {\bf w}^T_2 & 0 & 0& 0 & 0\\
        {\bf w}^T_3 & 0 & 0 & 0 & 0  
    \end{pmatrix}$$
    {\bf Claim}: $n(A(G))=3$.
    
    Here, \begin{center}
        
   ${\bf w}_1=\begin{pmatrix}
        {\bf 1}_{4t}\\
        {\bf 0}_{4s}\\
        {\bf 0}_{4\ell}
    \end{pmatrix}$, ${\bf w}_2=\begin{pmatrix}
        {\bf 0}_{4t}\\
        {\bf 1}_{4s}\\
        {\bf 0}_{4\ell}
    \end{pmatrix}$, ${\bf w}_3=\begin{pmatrix}
        {\bf 0}_{4t}\\
        {\bf 0}_{4s}\\
        {\bf 1}_{4\ell}
    \end{pmatrix}$.
     \end{center}
    where, $|N(u_{4k+2})|=4t,~|N(u_{4k+3})|=4s$ and $|N(u_{4k+2})|=4\ell,$ for some $t,s,\ell$ and $k=t+s+\ell.$\\
    Recall,  $x=(x_i)=(1,0,-1,0,1,0,\ldots,-1,0)^T,~y=(y_i)=(0,1,0,-1,0,\ldots,0,-1)^T$ are in $N(A(C_{4k}))$. We will use vectors $x,y$ to obtain null vectors for $A(G).$
\begin{eqnarray*}
    A(G)\begin{pmatrix}
        x\\
        0
    \end{pmatrix}&  = & \begin{pmatrix}
        A(C_{4k}) & {\bf 1} & {\bf w}_1 & {\bf w}_2 & {\bf w}_3\\
       { \bf 1}^T & 0 & 0 & 0 & 0\\
        {\bf w}^T_1 & 0 & 0 & 0 & 0\\
        {\bf w}^T_2 & 0 & 0& 0 & 0\\
        {\bf w}^T_3 & 0 & 0 & 0 & 0  
    \end{pmatrix}\begin{pmatrix}
        x_1\\
        \vdots\\
        x_{4k}\\
        {\bf 0}_{4}
    \end{pmatrix}\\
    & = & 0.
\end{eqnarray*} This is obtained from the following conclusions:\\
Using Lemma \ref{cyclenull}, we obtain $$A(C_{4k})x=0.$$ We can also observe that for any $i$, we have the following  $$\sum_{j=0}^{3}x_{i+j}=0.$$
Thus, we conclude that $w_i^Tx=0$ for each $1\le i \le 3$ and ${\bf 1}^Tx=0.$

On similar lines, one can conclude that $\begin{pmatrix}
    y\\0
\end{pmatrix}$ also belongs to the null space of $A(G).$ 
It is clear that $\begin{pmatrix}
    x\\0
\end{pmatrix}, \begin{pmatrix}
    y\\0
\end{pmatrix}$ and $z$ are linearly independent set of vectors. Therefore, $n(A)\ge 3.$ 

It remains to prove that $rk(A)=n-3$. Without loss of generality, assume that $u_1,u_2\in N(w_1)$. Then
consider the submatrix $A(1,2,4k+1)$ obtained by deleting the rows and columns indexed by the set $\{1,2,4k+1\}$  of $A$. We have $$B=A(1,2,4k+1)=\begin{pmatrix}
        A(P_{4t-2}) & D_{4t-2\times 4s} & {\bf 0} & {\bf 1}_{4t-2} & {\bf 0} & {\bf 0}\\
        D^T_{4t-2\times 4s} & A(P_{4s}) & D_{4s\times 4\ell} &{\bf 0} & {\bf 1}_{4s} &{\bf 0}\\
        {\bf 0} & D^T_{4s\times 4\ell} & A(P_{4\ell})& {\bf 0} & {\bf 0} & {\bf 1}_{4\ell}\\
         {\bf 1}^T_{4t-2} & {\bf 0} &{\bf 0} & 0 & 0& 0\\
         {\bf 0} & {\bf 1}^T_{4s} &{\bf 0} & 0 & 0 & 0\\
         {\bf 0} & {\bf 0} & {\bf 1}^T_{4\ell} & 0 & 0 & 0
    \end{pmatrix}$$
    where, $D_{m\times n}=[e_m,0,\ldots,0]$, we can observe that the matrix \begin{center}
     $C=\begin{pmatrix}
        A(P_{4t-2}) & D_{4t-2\times 4s} & {\bf 0}\\
    D^T_{4t-2\times 4s} & A(P_{4s}) & D_{4s\times 4\ell} \\
        {\bf 0} & D^T_{4s\times 4\ell} & A(P_{4\ell})
    \end{pmatrix}$   
    \end{center} is the adjacency matrix of $P_{4t+4s+4\ell-2}$, which is invertible. Therefore, using the determinant formula for a block matrix in the definition \ref{detblockfor}, we obtain the following: \begin{center}
      $\det(B)=\det(C)\det\bigg(\begin{pmatrix}
        0 & 0 & 0\\
        0 & 0& 0
        \\0 & 0& 0
    \end{pmatrix}-\begin{pmatrix}
         {\bf 1}^T_{4t-2} & {\bf 0} & {\bf 0} \\
         {\bf 0}&{\bf 1}^T_{4s} & {\bf 0}\\
         {\bf 0} & {\bf 0} & {\bf 1}^T_{4\ell}
    \end{pmatrix}C^{-1}\begin{pmatrix}
        {\bf 1}_{4t-2} & {\bf 0} & {\bf 0}\\
        {\bf 0} & {\bf 1}_{4s} & {\bf 0}\\
        {\bf 0} & {\bf 0} & {\bf 1}_{4\ell}
    \end{pmatrix}\bigg).$  
    \end{center}
Then, using Lemma \ref{pathodd}, $h=2t+2s+2\ell-1$. Thus, $\frac{h-(2t-1)}{2}=\ell+s.$
\begin{eqnarray*}
 \begin{pmatrix}
        C^{-1}\begin{pmatrix}{\bf 1}_{4t-2}\\ {\bf 0}\\ {\bf 0} 
        \end{pmatrix}& C^{-1}\begin{pmatrix}
             {\bf 0} \\ {\bf 1}_{4s}\\ {\bf 0}
        \end{pmatrix}&C^{-1}\begin{pmatrix}
            {\bf 0} \\ {\bf 0} \\ {\bf 1}_{4\ell}
        \end{pmatrix}
    \end{pmatrix}=\begin{pmatrix}
        1 & 0 & 0\\
        1 & 0 & 0\\
        \beta_{t-1} & {\bf 0}_{4(t-1)} & {\bf 0}_{4(t-1)} \\
      \zeta_{s} & \alpha_{s} & {\bf 0}_{4s}\\
      \zeta_{\ell} & {\bf 0}_{4\ell} & \alpha_{\ell}
    \end{pmatrix}
    \end{eqnarray*}
    where, $\beta_s=(\beta,\ldots,\beta)^T$,~$\beta=(0,0,1,1)^T,~ \alpha=(0,1,0,1)^T$ and $ \zeta=(0,-1,0,1)^T.$
Therefore, 
\begin{eqnarray*}
\begin{pmatrix}
         {\bf 1}^T_{4t-2} & {\bf 0} & {\bf 0} \\
         {\bf 0}&{\bf 1}^T_{4s} & {\bf 0}\\
         {\bf 0} & {\bf 0} & {\bf 1}^T_{4\ell}
    \end{pmatrix}\begin{pmatrix}
        1 & 0 & 0\\
        1 & 0 & 0\\
        \beta_{t-1} & {\bf 0}_{4(t-1)} & {\bf 0}_{4(t-1)} \\
      \zeta_{s} & \alpha_{s} & {\bf 0}_{4s}\\
      \zeta_{\ell} & {\bf 0}_{4\ell} & \alpha_{\ell}
    \end{pmatrix}&=&
    \begin{pmatrix}
       2+\sum_{t-1}{\bf 1}^T_4\beta & 0 & 0\\
       *  & \sum_{s}{\bf 1}^T_4\alpha  & 0\\
        * & * & \sum_{\ell}{\bf 1}^T_4\alpha 
    \end{pmatrix}\\&=&\begin{pmatrix}
        2t&0& 0\\
        * & 2s & 0\\
        * & * & 2\ell
    \end{pmatrix}
    \end{eqnarray*}
Thus, $\det(A(G))=-8ts\ell\det(C)$. We have $4t+4s+4\ell-2\cong 2(mod~4)$, using equation \ref{detpath}, we obtain $\det(C)=-1.$ Therefore, $\det(B)=8ts\ell,$ which is nonzero implies that $B$ is invertible. Indeed, this would imply $rk(A)\ge n-3.$  
This concludes $n(A)=3.$ 
 
  Thus, for each $i\in V,$ at least the $i$-th coordinate of one vector from the basis is nonzero. Therefore, $G$ is a core graph.
  \end{proof}

\begin{thm}
    Let $G\in$ \( \mathscr{C}_1\). Then the following holds:
    \begin{enumerate}
        \item $G$ satisfies the problem \ref{conj1},
        \item We define a vector $s=(s_i)\in R(A)$, for a fixed $t_o,~u_{t_o}\in V_2$,  $$s_i=\begin{dcases}
            1, & u_i\in \{u_{4k+1},u_{t_o}\} \\
            0 & Otherwise
        \end{dcases}.$$ Moreover, $G$ satisfies the problem \ref{conj2}.
\end{enumerate}
\end{thm}
\begin{proof}
\begin{enumerate}
    \item 
In the above proof, we have obtained a basis that consists of $\{-1,0,1\}$-vectors. Thus, $G$ satisfies the problem \ref{conj1}.
\item We can trivially see that the vector $s$ is orthogonal to the basis vectors of the null space of $A(G)$. Thus, it follows that $s\in R(A).$ The vector $s$ is not a column of $A(G),$ since every vertex adjacent to $u_{4k+1}$ is of degree $4.$ This establishes the claim.
\end{enumerate}
\end{proof}

The following example is an illustration of Construction II.

\begin{ex}\label{wheel38}
Let us consider the graph in Figure \ref{fig:7}. The partition of the vertex set given by
   $$V_0=\{1,2,\ldots,8\},~V_1=\{9\},~V_2=\{10,11\}.$$ 
  The neighborhoods of vertices in $V_2$ are defined as
$$N(10)=\{5,6,7,8\},~N(11)=\{1,2,3,4\}.$$
It follows that the collection $\{N(v):v\in V_2\}$ forms a partition of $V_0.$

Thus, calculating $N(A)$ using MATLAB, we obtained the following:$$\{(1,0,-1,0,1,0,-1,0,0,0),(0,-1,0,1,0,-1,0,1,0,0,0),(0,0,0,0,0,0,0,0,1,-1,-1)\},$$ which coincides with the null vectors obtained in Theorem \ref{C1core}. 
\end{ex}

\begin{rem}
    Suppose $G\in \mathscr{C}_1$ then clearly $G\in \mathscr{C}.$ By Theorem \ref{C1core}, $G$ is a core graph. Moreover, from the proof of Theorem \ref{C1core}, there exists a vector $z$ in the null space of $A$ that is not orthogonal to the all-ones vector. Hence, $0$ is a main eigenvalue of $G$.
    Clearly, $G$ is not regular, since  $d(u_{4k+1})=4k$ while $d(u_i)=4$ for $1\le i \le 4k$. Furthermore, since the subgraph induced by $V_0\cup V_1$ is a wheel graph, $G$ contains an odd cycle of length $3$, and therefore $G$ is not bipartite.
    If $|V_2|\le 3$, then by Lemma \ref{C1diam} we have $diam(G)=3.$ From definition \ref{Corewheel1} we conclude that the conditions $(3)$ and $(4)$ of Theorem \ref{sciriha} hold trivially. Thus, $G\in \mathscr{C}.$
\end{rem}

The following family of graphs, which is introduced here, helps us to obtain construction II.

\begin{defn}
    Let $n$ be a positive integer. We define a family of graphs $\mathscr{F}$ whose vertex set is partitioned into two subsets, $V_0$ and $V_1$ given by 
 $$V_0=\{u_1,\ldots u_{4k}\},~ V_1=\{u_{4k+1},\ldots,u_{4k+4s}\},$$ for integers $k,s\ge 1.$
     The edge set $E$ of $\mathscr{F}$ is defined as follows: 
     \begin{enumerate}
         \item The vertices in $V_0$ induce a cycle. \item The vertices in $V_1$ induce a cycle.
         \item The vertices $u_1,u_2$ are adjacent to both $u_{4k+4s-1}$ and $u_{4k+4s}$.
     \end{enumerate}
\end{defn}

\begin{lem}\label{null2c}
    If $G\in \mathscr{F}$, then $G$ is a core graph.
\end{lem}

\begin{proof}
    We have $$A(G)=\begin{pmatrix}
        A(C_{4k}) & E\\
        E^T & A(C_{4s})
    \end{pmatrix}$$
    where, $E=\begin{pmatrix}
        0 & \ldots & 0 & 1 & 1\\
        0 & \ldots & 0 & 1 &1\\
        {\bf 0} & \ldots& {\bf 0}& {\bf 0}& {\bf 0}\\
    \end{pmatrix}.$
    
     From Lemma \ref{cyclenull}, we obtain $x,y\in N(A(C_{4k}))$ and $p,q\in N(A(C_{4s})). $ We will use $x,y,p$ and $q $ to obtain null vectors of $A(G).$ 
    
     {\bf Claim}: $\begin{pmatrix}
        x-y\\
        0
    \end{pmatrix}, \begin{pmatrix}
        0 \\
        p-q
    \end{pmatrix}$ form a basis for $ N(A).$
    
     We have $A(C_{4k})(x-y)=0$. 
    We also have,$$E^T(x-y)=\begin{pmatrix}
         0& 0& {\bf 0}\\\vdots & \vdots & \vdots \\0 & 0 & {\bf 0}\\
       1 & 1 & {\bf 0}\\
         1 & 1 & {\bf 0}\\
       
    \end{pmatrix}\begin{pmatrix}
        ~~1 \\
        -1 \\
        -1\\\vdots \\
        ~~1
    \end{pmatrix}=0.$$
    This implies $\begin{pmatrix}
        x-y\\
        0
    \end{pmatrix}\in N(A)$. A similar argument concludes that $\begin{pmatrix}
        0 \\
        p-q
    \end{pmatrix}\in N(A).$ 
    
     It suffices to prove that the nullity of $A$ is $2$. We have already established that $n(A)\ge 2,$ which implies $rk(A)\le 4k+4s-2$. \\
   It remains to prove that $rk(A)=4k+4s-2.$
    Consider the submatrix $A(1,4k+4s)$ obtained by deleting the first and last columns and rows of $A$. It is enough to show that this submatrix is invertible. Indeed, this would imply $rk(A)\ge 4k+4s-2.$ This establishes the claim. It is easy to see that $A(1,4k+4s)$ is the adjacency matrix of a path on $4k+4s-2$ vertices. Since $4k+4s-2$ is even, its adjacency matrix is nonsingular. Thus, $A(1,4k+4s)$ is invertible and therefore $rk(A)=4k+4s-2.$ 
Thus, on the basis of $N(A)$ that we obtained, we conclude that $G$ is a core graph.
\end{proof}

\begin{rem}
    From the above lemma, we can observe that for any $G\in \mathscr{F},$ we have $$\begin{pmatrix}
        x-y\\
        0
    \end{pmatrix}, \begin{pmatrix}
        0 \\
        p-q
    \end{pmatrix}\in N(A),$$ and form a basis, which are $\{-1,0,1\}$-vectors. Thus, satisfies the problem \ref{conj1}.
\end{rem}

The following example is an illustration of the above family of graphs.

\begin{ex}\label{wheel38}
Consider the following graph \begin{figure}[H]
\begin{center}
\begin{tikzpicture}[scale=1]
\draw  (6.,5.)-- (7.,4.);
\draw  (7.,2.)-- (6.,0.8);
\draw  (4.,1.)-- (3.,2.);
\draw  (3.,4.)-- (4.,5.);
\draw  (6.,5.) -- (4.,5.);
\draw  (6.,0.8)--(4.,1.);
 \draw  (3.,4.)--(3.,2.);
 \draw (7.,4.)--(7.,2.) ;
\begin{scriptsize}
\fill (6.,5.) circle (2.5pt);
\draw (6.,5.25) node {$8$};
\fill (7.,4.) circle (2.5pt);
\draw (7.,4.3) node {$1$};
\fill (7.,2.) circle (2.5pt);
\draw (7.15,1.7) node {$2$};
\fill (6.,0.8) circle (2.5pt);
\draw (6.25,0.7) node {$3$};
\fill (4.,1.) circle (2.5pt);
\draw (3.75,0.7) node {$4$};
\fill (3.,2.) circle (2.5pt);
\draw (2.6,2.) node {$5$};
\fill (4.,5.) circle (2.5pt);
\draw (3.9,5.3) node {$7$};
\fill (3.,4.) circle (2.5pt);
\draw (2.7,4.) node {$6$};
\fill (11.,4.) circle (2.5pt);
\draw (11.,4.3) node {$9$};
\fill (11.,2.) circle (2.5pt);
\draw (11.,1.75) node {$10$};
\fill (9.,2.) circle (2.5pt);
\draw (9.,1.75) node {$11$};
\fill (9.,4.) circle (2.5pt);
\draw (9.,4.3) node {$12$};
\draw (9.,4.)--(11.,4.);
\draw (11.,4.)--(11.,2.);
\draw (11.,2.)--(9.,2.);
\draw (9.,2.)--(9.,4.);
\draw (9.,4.)--(7.,4.);
\draw (9.,2.)--(7.,2.);
\draw (9.,2.)--(7.,4.);
\draw (9.,4.)--(7.,2.);
\end{scriptsize}
\end{tikzpicture}
\end{center}
\caption{}
\label{fig:2wheels}
\end{figure} 
We have $$V_0=\{1,2,\ldots,8\},~V_1=\{9,10,11,12\}.$$ 
Thus, calculating $N(A)$ using MATLAB, we obtained the following:$$\{(1,-1,-1,1,1,-1,-1,1,0,0,0,0),(0,0,0,0,0,0,0,0,1,-1,-1,1)\}.$$ Therefore, these null vectors coincide with the null vectors obtained in Theorem \ref{null2c}, which verifies the result. 

\end{ex}

\subsection{Construction II}
Using the family of graphs $\mathscr{F}$, we obtain the following construction of core graphs.

\begin{defn}\label{W2}
Let $n$ be a positive integer. We define a family of graphs  $\mathscr{C}_2$ whose vertex set $V$ is partitioned into five subsets $V_0, V_1,V_2,V_3$ and $V_4,$ given by 
$$ V_0=\{u_1,\ldots, u_{4k}\},~V_1=\{u_{4k+1},\ldots, u_{4k+4s}\},$$
$$V_2=\{u_{4k+4s+1}\},~V_3=\{u_{4k+4s+2}\},~V_4=\{u_{4k+4s+3}\}.$$
    for integers $k,s\ge 1.$ The edge set \( E \) of \( \mathscr{C}_2\) is defined as follows:
\begin{enumerate}
\item The subgraph induced by $V_0\cup V_1$ is a graph  $G\in\mathscr{F}.$
\item Each vertex in $V_0$ is adjacent to $u_{4k+4s+1}\in V_2,$ and each vertex in $V_1$ is adjacent to $u_{4k+4s+2}\in V_3.$  
    \item The vertex $w\in V_4$ is adjacent to every vertex in $V_0\cup V_1$.
\end{enumerate}
\end{defn}

\begin{figure}[H]
\begin{center}
\begin{tikzpicture}[scale=1]
\draw  (4.7,3.)-- (6.,5.);
\draw  (6.,5.)-- (7.,4.);
\draw  (7.,4.)-- (4.7,3.);
\draw  (4.7,3.)-- (7.3,2.);
\draw  (7.3,2.)-- (5.7,1.);
\draw  (5.7,1.)-- (4.7,3.);
\draw  (4.7,3.)-- (4.3,1.);
\draw  (4.3,1.)-- (3.,2.);
\draw  (3.,2.)-- (4.7,3.);
\draw  (4.7,3.)-- (3.,3.5);
\draw  (3.,3.5)-- (4.,5.);
\draw  (4.,5.)-- (4.7,3.);
\draw  (6.,5.) -- (4.,5.);
\draw  (5.7,1.)--(4.3,1.);
 \draw  (3.,3.5)--(3.,2.);
 \draw (7.,4.)--(7.3,2.) ;
 \draw (8.,8.)--(6.,5.);
 \draw (8.,8.)--(3.,2.);
 \draw (8.,8.)--(3.,3.5);
 \draw (8.,8.)--(4.3,1.);
 \draw (8.,8.)--(4.,5.);
 \draw (8.,8.)--(5.7,1.);
 \draw (8.,8.)--(6.,5.);
\begin{scriptsize}
\fill (8.,8.) circle (2.5pt);
\draw (8.,8.35) node {$19$};
\fill (4.7,3.) circle (2.5pt);
\draw (4.75,3.45) node {$17$};
\fill (6.,5.) circle (2.5pt);
\draw (5.82,5.2) node {$8$};
\fill (7.,4.) circle (2.5pt);
\draw (6.9,3.7)
node {$1$};
\fill (7.3,2.) circle (2.5pt);
\draw (7.3,1.7) node {$2$};
\fill (5.7,1.) circle (2.5pt);
\draw (5.7,0.7) node {$3$};
\fill (4.3,1.) circle (2.5pt);
\draw (4.2,0.7) node {$4$};
\fill (3.,2.) circle (2.5pt);
\draw (2.75,2.) 
node {$5$};
\fill (4.,5.) circle (2.5pt);
\draw (3.9,5.3) node {$7$};
\fill (3.,3.5) circle (2.5pt);
\draw (2.7,3.6) 
node {$6$};
\fill (11.3,3.) circle (2.5pt);
\draw (11.25,3.45) node {$18$};
\fill (12.,5.) circle (2.5pt);
\draw (12.25,5.25) node {$10$};
\fill (13.,3.5) circle (2.5pt);
\draw (13.35,3.65) node {$11$};
\fill (13.,2.) circle (2.5pt);
\draw (13.2,1.8) node {$12$};
\fill (11.6,1.) circle (2.5pt);
\draw (11.7,0.7) node {$13$};
\fill (10.4,1.) circle (2.5pt);
\draw (10.3,0.7) node {$14$};
\fill (8.7,2.) circle (2.5pt);
\draw (8.6,1.75) node {$15$};
\fill (10.,5.) circle (2.5pt);
\draw (9.7,5.1) node {$9$};
\fill (9.,4.) circle (2.5pt);
\draw (8.7,4.3) node {$16$};
\draw (8.,8.)--(10.,5.);
\draw (8.,8.)--(10.4,1.);
\draw (8.,8.)--(8.7,2.);
\draw (8.,8.)--(11.6,1.);
\draw (8.,8.)--(13.,2.);
\draw (8.,8.)--(13.,3.5);\draw (8.,8.)--(9.,4.);
\draw (8.,8.)--(12.,5.);
\draw (8.,8.)--(7.,4.);
\draw (8.,8.)--(7.3,2.);
\draw  (11.3,3.)--(10.,5.);
\draw  (11.3,3.)--(8.7,2.);
\draw  (11.3,3.)--(10.4,1.);
\draw  (11.3,3.)--(11.6,1.);\draw  (11.3,3.)--(13.,2.);
\draw  (11.3,3.)--(13.,3.5);
\draw  (11.3,3.)--(10.,5.);
\draw (11.3,3.)--(12.,5.);
\draw (11.3,3.)--(9.,4.);
\draw (12.,5.)--(10.,5.) ;
\draw (12.,5.)--(13.,3.5);
\draw (13.,3.5)--(13.,2.);
\draw (10.,5.)--(9.,4.);
 \draw (13.,2.)--(11.6,1.);
 \draw (11.6,1.)--(10.4,1.);
  \draw (10.4,1.)--(8.7,2.);
   \draw (8.7,2.)--(9.,4.);
   \draw (9.,4.)--(7.,4.);
    \draw (8.7,2.)--(7.3,2.);
     \draw (8.7,2.)--(7.,4.);
      \draw (9.,4.)--(7.3,2.);
\end{scriptsize}
\end{tikzpicture}
\end{center}
\caption{}
\label{fig:2wheels}
\end{figure}

\begin{lem}\label{c2diam}
If $G\in$ \( \mathscr{C}_2\) then $diam(G)=3.$
\end{lem}

\begin{proof} Let $x,y\in V$. Assume $x,y$ are non-adjacent. It is enough to prove $d(x,y)\le 3$.\\
Let us assume $w=u_{4k+4s+3}$.\\
{\bf Case~1}: For $x,y\in V_0\cup V_1.$ \\
From the definition \ref{W2}~$(3)$, we conclude that $x,y\in N(w).$ Therefore, $d(x,y)=2.$ \\
{\bf Case~2}: For $x=u_{4k+4s+1},~y\in V_1.$\\
Here, there exists a path from $u_{4k+4s+1} $ to $y,$ passing through $u_i\in V_0$ for any $i$ and  vertex $w$, i.e.,  
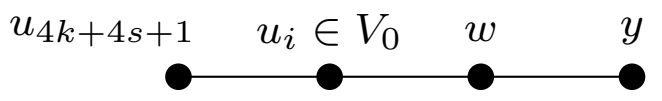
\begin{figure}[H]
\begin{center}
\begin{tikzpicture}[scale=1]
\draw (5.,3.)--(6.,3.);
   \draw (7.,3.)--(6.,3.);
   \draw (7.,3.)--(8.,3.);
   
\begin{scriptsize}
\fill (5.,3.) circle (2.5pt);
\draw (4.5,3.3) node {$u_{4k+4s+1}$};
\fill (6.,3.) circle (2.5pt);
\draw (6.,3.3) node {$u_i\in V_0$};
\fill (7.,3.) circle (2.5pt);
\draw (7.,3.3) node {$w$};
\fill (8.,3.) circle (2.5pt);
\draw (8.,3.3) node {$y$};

\end{scriptsize}
\end{tikzpicture}
\end{center}

\end{figure}
Therefore, $d(u_{4k+4s+1},y)\le 3.$\\ 
{\bf Case~3}: For $x=u_{4k+4s+2},~y\in V_0.$\\Here, there exists a path from $u_{4k+4s+2}$ to $y$ passing through $u_i\in V_1$ for any $i$ and  vertex $w$, i.e., 
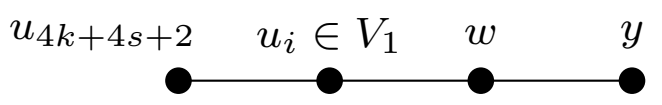
\begin{figure}[H]
\begin{center}
\begin{tikzpicture}[scale=1]
\draw (5.,3.)--(6.,3.);
   \draw (7.,3.)--(6.,3.);
   \draw (7.,3.)--(8.,3.);
   
\begin{scriptsize}
\fill (5.,3.) circle (2.5pt);
\draw (4.5,3.3) node {$u_{4k+4s+2}$};
\fill (6.,3.) circle (2.5pt);
\draw (6.,3.3) node {$u_i\in V_1$};
\fill (7.,3.) circle (2.5pt);
\draw (7.,3.3) node {$w$};
\fill (8.,3.) circle (2.5pt);
\draw (8.,3.3) node {$y$};

\end{scriptsize}
\end{tikzpicture}
\end{center}
\end{figure}
 Thus, $d(u_{4k+4s+2},y)\le 3.$\\
{\bf Case~4}: For $x=u_{4k+4s+1},~y=u_{4k+4s+2}$.\\
We can observe that there exists a path from $u_{4k+4s+1}$ to $u_{4k+4s+2},$ passing through the vertices $u_1,u_{4k+4s-1}$, which is of length $3.$ So that, $d(u_{4k+4s+1},u_{4k+4s+2})\le 3.$\\
{\bf Case ~5}: For $x=w,~y\in \{u_{4k+4s+1},u_{4k+4s+2}\}$.\\
Suppose $y=u_{4k+4s+1}$, then there exists a path from $w$ to $u_{4k+4s+1}$, passing through $u_i\in V_0$ for any $i$. i.e., 
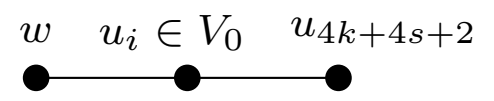
\begin{figure}[H]
\begin{center}
\begin{tikzpicture}[scale=1]
\draw (5.,3.)--(6.,3.);
   \draw (7.,3.)--(6.,3.);

\begin{scriptsize}
\fill (5.,3.) circle (2.5pt);
\draw (5.,3.3) node {$w$};
\fill (6.,3.) circle (2.5pt);
\draw (5.9,3.3) node {$u_i\in V_0$};
\fill (7.,3.) circle (2.5pt);
\draw (7.3,3.3) node {$u_{4k+4s+2}$};

\end{scriptsize}
\end{tikzpicture}
\end{center}
\end{figure}

Thus, we obtain $d(w,u_{4k+4s+1})=2.$ On the otherhand, if $y=u_{4k+4s+2},$ then there exists a path from $w$ to $u_{4k+4s+2}$, passing through $u_i\in V_1$ for any $i$. Thus, we obtain $d(w,u_{4k+4s+2})=2.$
Therefore, $d(x,y)\le 3$ for $x,y\in V,$ which establishes the claim.
\end{proof}

The following lemma is used to prove Theorem \ref{c2core}.

\begin{lem}\label{pathodddeg}
    Let $A$ be the adjacency matrix of the path graph $P_{2h}.$ Then the following statements hold:\begin{enumerate} \item If $A(P_{2h})x=\begin{pmatrix}1\\1\\ 1\\~~{\bf 0}_{m} \end{pmatrix}$
    then $x=\begin{pmatrix}
        b\\
        ~~{\bf 0}_m
    \end{pmatrix}$ where, $b=(1,1,0)^T.$
      \item If 
$A(P_{2h})x=\begin{pmatrix}~~{\bf 1}_{4t-1} \\ {\bf 0}_{m} \end{pmatrix}$
    then $x=\begin{pmatrix}
         b\\
        ~~\alpha_{t-1}\\
        ~~{\bf 0}_m
    \end{pmatrix}$ 
\item If
$A(P_{2h})x=\begin{pmatrix}
    {\bf 0}_m \\
    ~~{\bf 1}_{4s-1}  \end{pmatrix}$ then $x=\begin{pmatrix}
        {\bf 0}_m\\
        ~~\alpha_{s-1}\\
        g
    \end{pmatrix}$  where, $g=(0,1,1)^T$ and  $\alpha_{i}$ is defined as in Lemma \ref{lem1}. 
    \end{enumerate}
\end{lem}

\begin{proof}
\begin{enumerate}
\item It is trivial to observe that $$A\begin{pmatrix}
    b\\
    ~~{\bf 0}_m\\
    \end{pmatrix}=A\begin{pmatrix}
        1\\
        1\\0\\
        ~~{\bf 0}_{m}
    \end{pmatrix}=Ae^{1}+Ae^{2}=\begin{pmatrix}1\\1\\ 1\\~{\bf 0}_{m} \end{pmatrix}$$
        \item

        We have 
        \begin{eqnarray*}
        A\begin{pmatrix}
        b\\~~\alpha_{t-1}\\
        ~~{\bf 0}_m
    \end{pmatrix}&=&A\begin{pmatrix}
        b\\
        ~~{\bf 0}_{4t-4}\\
        ~~{\bf 0}_{m}
    \end{pmatrix}+A\begin{pmatrix}
        {\bf 0}_3\\
        ~~\alpha_{t-1}\\
        ~~{\bf 0}_m
    \end{pmatrix} \\
    &=& \begin{pmatrix}
        {\bf 1}_3\\
        {\bf 0}_{4t-4}\\
        {\bf 0}_{m}
    \end{pmatrix}+\begin{pmatrix}
        {\bf 0}_3\\
        {\bf 1}_{4t-4}\\
        {\bf 0}_m
    \end{pmatrix}\\
    &=&\begin{pmatrix}
        ~~{\bf 1}_{4t-1}\\
        {\bf 0}_m
    \end{pmatrix}
        \end{eqnarray*}
      
    \item A similar argument establishes the claim. 
\end{enumerate}
\end{proof}

\begin{thm}\label{c2core}
Every graph $G\in$ \( \mathscr{C}_2\) is a core graph.
\end{thm}
\begin{proof}
 We have $$A(G)=\begin{pmatrix}
     A(C_{4k}) & E & {\bf 1} & {\bf 0} & {\bf 1} \\
     ~~E^T & A(C_{4s}) & {\bf 0} & {\bf 1} & {\bf 1} \\
     ~~{\bf 1}^T & {\bf 0} & 0 & 0 & 0\\
     {\bf 0} & ~~{\bf 1}^T & 0 & 0 & 0\\
     ~~{\bf 1}^T & ~~{\bf 1}^T & 0 & 0 & 0 
\end{pmatrix}.$$
{\bf Claim}:  $\begin{pmatrix}
    x-y  \\ {\bf 0} \\0 \\0 \\0
\end{pmatrix},\begin{pmatrix}
   {\bf 0} \\p-q \\0\\0\\0
\end{pmatrix},\begin{pmatrix}
    ~~{\bf 0}\\
    ~~{\bf 0}\\
    ~~1\\
    ~~1\\
    -1
\end{pmatrix}$ form a basis for $N(A).$

It is easy to verify that these three vectors belong to the null space of $A$ by following the arguments in Theorem \ref{C1core} and Lemma \ref{null2c}. Since these vectors are linearly independent, $n(A)\ge 3.$ This implies $rk(A)\le 4k+4s.$ It remains to prove that $rk(A)=4k+4s.$\\
    Consider the submatrix $B=A(1,4k+4s,4k+4s+3)$ obtained by deleting the first, $4k+4s$ and $4k+4s+3$ rows and columns of $A$. 
    We have, $$B=\begin{pmatrix}
        A(P_{4t-1}) & D& {\bf 1}_{4t-1} & {\bf 0}\\
        D^T & A(P_{4s-1}) & {\bf 0} & {\bf 1}_{4s-1}\\
         {\bf 1}^T_{4t-1} & {\bf 0} & 0 & 0\\
         {\bf 0} & {\bf 1}^T_{4s-1} & 0 & 0
    \end{pmatrix}$$
    where, $D^T=[e_{4s-1},0,\ldots,0]$, we can observe that the matrix $$C=\begin{pmatrix}
        A(P_{4t-1}) & D\\
        D^T & A(P_{4s-1})
    \end{pmatrix}$$ is the adjacency matrix of $P_{4t+4s-2}$, which is invertible. Therefore, using the determinant formula for a block matrix, we obtain the following $$\det(A(G))=\det(C)\det(\begin{pmatrix}
        0 & 0\\
        0 & 0
    \end{pmatrix}-\begin{pmatrix}
         {\bf 1}^T_{4t-1} & {\bf 0} \\
         {\bf 0}&{\bf 1}^T_{4s-1} 
    \end{pmatrix}C^{-1}\begin{pmatrix}
        {\bf 1}_{4t-1} & {\bf 0}\\
        {\bf 0} & {\bf 1}_{4s-1}
    \end{pmatrix}).$$
    
    Here, $b=(1,1,0)^T$ and $g=(0,1,1)^T.$ Using Lemma \ref{pathodddeg}, we have the following:
    \begin{eqnarray*}
    \begin{pmatrix}
         {\bf 1}^T_{4t-1} & {\bf 0} \\
         {\bf 0}&{\bf 1}^T_{4s-1} 
    \end{pmatrix}C^{-1}\begin{pmatrix}
        {\bf 1}_{4t-1} & {\bf 0}\\
        {\bf 0} & {\bf 1}_{4s-1}
    \end{pmatrix} &= & \begin{pmatrix}
         {\bf 1}^T_{4t-1} & {\bf 0} \\
         {\bf 0}&{\bf 1}^T_{4s-1} 
    \end{pmatrix}\begin{pmatrix}
C^{-1}\begin{pmatrix}
           {\bf 1}_{4t-1}\\ {\bf 0}_{4s-1}  
\end{pmatrix}&
C^{-1}\begin{pmatrix}  {\bf 0}_{4t-1}\\{\bf 1}_{4s-1} 
    \end{pmatrix}
    \end{pmatrix}\\
&=&\begin{pmatrix}
         {\bf 1}^T_{4t-1} & {\bf 0}_{4s-1} \\
         {\bf 0}_{4t-1}&{\bf 1}^T_{4s-1} 
    \end{pmatrix}
        \begin{pmatrix}
        b & {\bf 0}_3\\
        \alpha_{t-1}&{\bf 0}_{4t-4} \\
        {\bf 0}_{4s-4} & \alpha_{s-1}\\{\bf 0}_3& g
    \end{pmatrix}\\
    & = &\begin{pmatrix}
        {\bf 1}_3^Tb+\sum_{t-1}({\bf 1}_4^T\alpha) &0\\0&\sum_{s-1}({\bf 1}_4^T\alpha)+{\bf 1}_3^Tg
    \end{pmatrix}\\
    &=&\begin{pmatrix}
        2t & 0\\0& 2s
    \end{pmatrix}\end{eqnarray*}

   Therefore, from the above arguments, we conclude that
   \begin{eqnarray*}
   \det(B)&=&\det(C)\det\begin{pmatrix}
       -2t & 0\\0 & -2s
   \end{pmatrix}\\
   &=&4ts\det(C).
   \end{eqnarray*} We have $4t+4s-2\cong 2(mod~4)$, using equation \ref{detpath}, we obtain $\det(C)=-1.$
    Thus, $\det(B)=-4ts$, nonzero, implies that $B$ is invertible. Indeed, this would imply $rk(A)\ge 4k+4s.$ Thus, the nullity is $3$. This establishes the claim. 

   Thus, for each $i\in V,$ at least the $i$-th coordinate of one vector from the basis is nonzero. Therefore, $G$ is a core graph.
\end{proof}


\begin{thm}\label{nullc2graph}
    Let $G\in$ \( \mathscr{C}_2\). Then the following holds:
    \begin{enumerate}
        \item $G$ satisfies the problem \ref{conj1},
        \item  We define a vector $s=(s_i)\in R(A)$, $$s_i=\begin{dcases}
            1, & i=4k+4s+1,4k+4s+3, \\
            0 & otherwise.
        \end{dcases}$$ Moreover, $G$ satisfies the problem \ref{conj2}.
\end{enumerate}
\end{thm}
\begin{proof}
\begin{enumerate}
    \item
    We can observe that $\begin{pmatrix}
    x-y  \\
    {\bf 0} \\
    0 \\
    0 \\
    0
\end{pmatrix},\begin{pmatrix}
   {\bf 0}\\
   p-q\\
   0\\
   0\\
   0
\end{pmatrix},\begin{pmatrix}
    ~~{\bf 0}\\
    ~~{\bf 0}\\
    ~~1\\
    ~~1\\
    -1
\end{pmatrix}$ form the basis for $N(A),$ which are $\{-1,0,1\}$-vectors. Thus, $G$ satisfies the problem \ref{conj1}.
\item We can easily observe that the vector $s$ is orthogonal to the basis vectors. Thus, $s\in R(A),$ which is a $(0,1)$-vector. Thus, $G$ satisfies the problem \ref{conj2}.
\end{enumerate}
\end{proof}

\begin{rem}
    Suppose $G\in \mathscr{C}_2$  then $G\in \mathscr{C}.$ By Theorem \ref{c2core}, $G$ is a core graph. Moreover, from the proof of Theorem \ref{c2core}, there exists a vector
    $({\bf 0},{\bf 0},1,1,-1)^T$ in the null space of $A$ that is not orthogonal to the all-ones vector. Hence, $0$ is a main eigenvalue of $G$. Clearly, $G$ is not regular, since  $d(u_{4k+4s+1})=4k$, while $d(u_i)=4,$ for $3\le i \le 4k$. Furthermore, since the subgraph induced by $V_0\cup V_2$ is a wheel graph, $G$ contains an odd cycle of length $3$, and therefore $G$ is not bipartite. By Lemma \ref{c2diam} we have $diam(G)=3.$ 
    From the definition \ref{W2} we conclude that conditions $(3)$ and $(4)$ of Theorem \ref{sciriha} are trivially valid. Thus, $G\in \mathscr{C}.$
\end{rem}

The following example is an illustration of construction II.

\begin{ex}\label{wheel38}
Consider the graph in figure  \ref{fig:2wheels}. We have $$V_0=\{1,2,\ldots,8\},~V_1=\{9,10,\ldots,16\},~V_2=\{17\},$$
$$V_3=\{18\},V_4=\{19\}.$$ We can observe that the subgraph induced by vertex set $V_0 \cup V_1=\{1,2\ldots,16\},$ belongs to $\mathscr{F}.$ Other conditions $(2),(3)$ from definition \ref{W2} holds trivally. 
Thus, calculating $N(A)$ using MATLAB, we obtain the following:$$\{(1,-1,-1,1,1,-1-1,1,{\bf 0}_{11}),({\bf 0}_8,1,-1,-1,1,1,-1,-1,1,{\bf 0}_3),({\bf 0}_{16},1,1,-1)\}.$$ Therefore, these null vectors coincide with the null vectors obtained in Theorem \ref{null2c}, which verifies the result. 
\end{ex}

\section{Some nut graphs satisfying the ACK conjecture}

The following result states that every nut graph with diameter $3$ satisfies the problem \ref{conj2}. 

\begin{thm}
    Let $G$ be a nut graph on $n$ vertices. If $diam(G)=3$, then $G$ satisfies the problem \ref{conj2}.
\end{thm}

\begin{proof}
    We have $diam(G)=3$, so there exists $u,v\in V(G)$ such that $d(u,v)=3.$ Since $u,v$ are not adjacent to a common vertex, $a^u+a^v$ is a (0,1) vector. \\
    {\bf Claim}: $a^u+a^v$ is not a column of $A.$ \\
     Suppose $a^u+a^v=a^\ell$ for some $\ell\in V(G)$ then $y:=e^u+e^v-e^\ell\in N(A).$ Since $diam(G)=3$, we have $n\ge 4$. Because $n(A)=1,$ it follows that every null vector of $G$ is of the form $\alpha y$ for some $\alpha\neq 0.$ Therefore, every null vector of $G$ has exactly three nonzero coordinates, contradicting the fact that $G$ is a core graph. \end{proof}

\section*{Acknowledgements}

The authors sincerely thank S. Akansha and Aditi Howlader for their valuable discussions and thoughtful insights that helped shape this work.

\end{document}